\documentclass[12pt]{article}

\usepackage{amsmath,amssymb}
\usepackage{color}
\newtheorem{theorem}{Theorem}[section]

\newtheorem{proposition}[theorem]{Proposition}

\newcommand{\proof}{\noindent{\bf Proof.\ }}
\newcommand{\qed}{\hfill $\square$ \bigskip}
\def\cp{\,\square\,}
\DeclareMathOperator {\gp} {gp}
\DeclareMathOperator {\diam} {diam}

\begin{document}

\title{The general position problem on Kneser graphs and on some graph operations}

\author{Modjtaba Ghorbani $^{a}$
\and
Sandi Klav\v zar $^{b,c,d}$
\and
Hamid Reza Maimani $^{a}$ 
\and 
Mostafa Momeni $^{a}$ 
\and 
Farhad Rahimi Mahid $^{a}$
\and
Gregor Rus $^{e}$}

\date{\today}

\maketitle

\begin{center}
$^{a}$ Department of Mathematics, Faculty of Science, Shahid Rajaee \\ 
Teacher Training University, Tehran, 16785 - 136, I. R. Iran
\medskip

$^b$ Faculty of Mathematics and Physics, University of Ljubljana, Slovenia\\
\medskip

$^c$ Faculty of Natural Sciences and Mathematics, University of Maribor, Slovenia\\
\medskip

$^d$ Institute of Mathematics, Physics and Mechanics, Ljubljana, Slovenia\\
\medskip

$^e$ Faculty of Organizational Sciences, University of Maribor, Slovenia\\
\end{center}

\begin{abstract}
A vertex subset $S$ of a graph $G$ is a general position set of $G$ if no vertex of $S$ lies on a geodesic between two other vertices of $S$. The cardinality of a largest general position set of $G$ is the general position number (gp-number) ${\rm gp}(G)$ of $G$. The gp-number is determined for some families of Kneser graphs, in particular for $K(n,2)$ and $K(n,3)$. A sharp lower bound on the gp-number is proved for Cartesian products of graphs. The gp-number is also determined for joins of graphs, coronas over graphs, and line graphs of complete graphs. 
\end{abstract}


\medskip\noindent
{\bf AMS Subj.\ Class:} 05C12; 05C69; 05C76 

\section{Introduction}
\label{sec:intro}

A general position problem in graph theory is to find a largest set of vertices that are in a general position. More precisely, if $G=(V(G), E(G))$ is a graph, then $S\subseteq V(G)$ is a {\em general position set} if for any triple of pairwise different vertices $u,v,w\in S$ we have $d_G(u,v) \ne d_G(u,w) + d_G(w,u$), where $d_G$ is the standard shortest path distance function in the graph $G$.  $S$ is called a {\em gp-set} of $G$ if $S$ has the largest cardinality among the general position sets of $G$. The {\em general position number} ({\em gp-number} for short) $\gp(G)$ of $G$ is the cardinality of a gp-set of $G$. 

This concept was introduced---under the present name---in~\cite{manuel-2018a} in part motivated by the  Dudeney's 1917 no-three-in-line problem ~\cite{dudeney-1917} (see~\cite{ku-2018, misiak-2016, PoWo07} for recent related results) and by a corresponding problem in discrete geometry known as the general position subset selection problem~\cite{froese-2017, payne-2013}. Independently geodetic irredundant sets were earlier introduced in~\cite{ullas-2016}, a concept which is equivalent to the general position sets. 

In~\cite{ullas-2016} graphs $G$ with $\gp(G)\in \{2, n(G)-1, n(G)\}$ were classified and some other results presented. Here and later, $n(G)$ denotes the order of $G$. Then, in~\cite{manuel-2018a}, several general bounds on the gp-number were presented, proved that set of simplicial vertices of a block graph form its gp-set, and proved that the problem is NP-complete in general. The gp-number of a large class of subgraphs of the infinite grid graph and of the infinite diagonal grid has been determined in~\cite{manuel-2018b}. In the paper~\cite{anand-2019+} a formula for the gp-number of graphs of diameter $2$ was given which in particular implies that $\gp(G)$ of a cograph $G$ can be determined in polynomial time. Moreover, a formula for the gp-number of the complement of a bipartite graph was also deduced. The main result of~\cite{anand-2019+} gives a characterization of general position sets (see Theorem~\ref{thm:gpsets} below). 

We proceed as follows. In the rest of this section further definitions are given and know results needed are stated. In Section~\ref{sec:Kneser} the gp-number is determined for some families of Kneser graphs. In particular, if $n\ge 7$, then $\gp(K(n,2)) = n-1$ and if $n\ge 8$, then $\gp(K(n,3)) = \binom{n-1}{2}$. In the subsequent section the gp-number of Cartesian products is bounded from below. The bound is proved to be sharp on the Cartesian product of complete graphs. We conclude the paper with Section~\ref{sec:more-opearations} in which the gp-number is determined for joins of graphs, coronas over graphs, and line graphs of complete graphs, where the first two results are stated as functions of the corresponding invariants of factor graphs. 

For a positive integer $n$ let $[n] = \{1,\ldots,n\}$. Graphs in this paper are finite, undirected, and simple. The maximum distance between all pairs of vertices of $G$ is the {\em diameter} $\diam(G)$ of $G$. An $u,v$-path of length $d_G(u,v)$ is called an $u,v$-{\it geodesic}. The {\em interval} $I_G[u,v]$ between vertices $u$ and $v$ of a graph $G$ is the set of vertices $x$ such that there exists a $u,v$-geodesic which contains $x$. A subgraph $H$ of $G$ is convex if for every $u,v\in V(H)$, all the vertices from $I_G(u,v)$ belong to $V(H)$. 

The size of a largest complete subgraph of a graph $G$ and the size of its largest independent set are denoted by $\omega(G)$ and $\alpha(G)$, respectively. The complement of a graph $G$ will be denoted with $\overline{G}$ and the subgraph of $G$ induced by $S\subseteq V(G)$ with $G[S]$. Let $\eta(G)$ denote the maximum order of an induced complete multipartite subgraph of $\overline{G}$. We will use the following result. 

\begin{theorem} {\rm \cite[Theorem 4.1]{anand-2019+}}
\label{thm:diameter2}
If $\diam(G) = 2$, then $\gp(G) = \max\{\omega(G), \eta(G)\}$. 
\end{theorem}

To complete the introduction we recall a characterization of general position sets from~\cite{anand-2019+}, for which some preparation is required. if $G$ is a connected graph, $S\subseteq V(G)$, and ${\cal P} = \{S_1, \ldots, S_p\}$ a partition of $S$, then ${\cal P}$ is \emph{distance-constant} (named ``distance-regular'' in \cite[p.~331]{kante-2017}) if for any $i,j\in [p]$, $i\ne j$, the distance $d_G(u,v)$, where $u\in S_i$ and $v\in S_j$, is independent of the selection of $u$ and $v$. This distance is then the distance $d_G(S_i,S_j)$ between the parts $S_i$ and $S_j$. A distance-constant partition ${\cal P}$ is {\em in-transitive} if $d_G(S_i, S_k) \ne d_G(S_i, S_j) + d_G(S_j,S_k)$ holds for arbitrary pairwise different $i,j,k\in [p]$. Then we have: 

\begin{theorem} {\rm \cite[Theorem 3.1]{anand-2019+}}
\label{thm:gpsets}
Let $G$ be a connected graph. Then $S\subseteq V(G)$ is a general position set if and only if the components of $G[S]$ are complete subgraphs, the vertices of which form an in-transitive, distance-constant partition of $S$. 
\end{theorem}

\section{Kneser graphs}
\label{sec:Kneser}

If $n$ and $k$ are positive integers with $n\ge k$, then the {\em Kneser graph} $K(n,k)$ has as vertices all the $k$-element subsets of the set $[n]$, vertices being adjacent if the corresponding sets are disjoint. For more on Kneser graph see~\cite{boruzan-2019, bresar-2019, mutze-2017, valencia-pabon-2005}. 

In this section we are interested in the gp-number of Kneser graphs, for which the following result will be useful.

\begin{proposition} {\rm \cite[Proposition 1]{valencia-pabon-2005}}
\label{prop:Kneser-diam-2}
If $k\ge 2$ and $n\ge 3k-1$, then $\diam(K(n,k)) = 2$.  
\end{proposition}

Recall also that the celebrated Erd\H{o}s-Ko-Rado theorem~\cite{erdos-1961} asserts that if $n\ge 2k$, then $\alpha(K(n,k)) \le \binom{n-1}{k-1}$, cf. also~\cite[Theorem 6.4]{van-lint-1992}. 

In our first result of the section we determine the gp-number of the Kneser graphs $K(n,2)$ as follows. 

\begin{theorem}  
\label{thm:kneser}
If $n\ge 4$, then 
\[\gp(K(n,2)) = \left\{ 
{\begin{array}{{ll}}
6; & 4\leq n\leq 6, \\
n-1; & n\geq 7\,.
\end{array}} \right.\]
\end{theorem}

\proof
Since $K(4,2)=3K_2$, we clearly have $\gp(K(4,2))=6$. The Kneser graph $K(5,2)$ is the Petersen graph for which it has been proven in~\cite{manuel-2018a} that $\gp(K(5,2)) = 6$. In the rest let $n\ge 6$. 

Let $S$ be a general position set of $K(n,2)$. By Theorem~\ref{thm:gpsets} the components of $K(n,2)[S]$ are complete graphs. Suppose one of this components, say $K$, has at least three vertices. Then any other vertex of $K(n,2)$ is adjacent to at least one vertex of $K$. It follows that $K(n,2)[S]$ has only one (complete) component and consequently $|n(K)| \le \lfloor\frac{n}{2}\rfloor$.

Suppose next that $K(n,2)[S]$ contains a component $K$ isomorphic to $K_2$. Assume without loss of generality that $V(K) = \{\{1,2\}, \{3,4\}\}$. Since no other vertex of $S$ is adjacent with the vertices of $K_2$, the other vertices of $S$ must be 2-subsets of $[4]$. It follows that if $K(n,2)[S]$ contains a component $K_2$, then 
$|S| \le 6$.

Assume finally that $K(n,2)[S]$ contains only isolated vertices, that is, $S$ is an independent set. Then the Erd\H{o}s-Ko-Rado theorem implies that $|S|\le n-1$. 

Suppose that $n = 6$. Then the six 2-subsets of $[4]$ induce three independent edges, hence $\gp(K(6,2)) \ge 6$. By the above we conclude that $\gp(K(6,2)) = 6$.

Let $n\ge 7$. Then by the above, $\gp(K(n,2)) \le n-1$. On the other hand, the set $\{\{1,2\}, \{1,3\}, \ldots, \{1,n\}\}$ is an independent set of $K(n,2)$. Since $\diam(K(n,2)) = 2$ by Proposition~\ref{prop:Kneser-diam-2}, we also have $\gp(K(n,2)) \ge n-1$. 
\qed

\begin{theorem}
\label{thm:Kneser-special}
Let $n, k \in \mathbb{N}$ and $n\ge 3k-1$. If for all $t$, where $2\le t \le k$, the inequality $k^t\binom{n-t}{k-t}+t\leq \binom{n-1}{k-1}$ holds, then 
$$\gp(K(n,k))=\binom{n-1}{k-1}\,.$$
\end{theorem} 


\proof
Since $n\ge 3k-1$, Proposition~\ref{prop:Kneser-diam-2} implies that $\diam(K(n,k)) = 2$.

Let ${\cal S}$ be the set of all $k$-subsets of $[n]$ that contain 1. Clearly, $|{\cal S}| = \binom{n-1}{k-1}$ and ${\cal S}$ form an independent set of $K(n,k)$. Hence, as  $\diam(K(n,k)) = 2$, we infer that ${\cal S}$ is a general position set and consequently $\gp(K(n,k)) \ge \binom{n-1}{k-1}$. 
 
Let $T$ be a general position set of $K(n,k)$ and let $K$ the a largest component of $K(n,k)[T]$. By Theorem~\ref{thm:gpsets} we know that $K$ is a complete subgraph. Let $|n(K)| = t$. If $t>k$, then every vertex $V(K(n,k))\setminus V(K)$ must have a neighbor in $K$. This implies that $T$ is the only component of $K(n,k)[T]$, but then we clearly have $|K| \le  \binom{n-1}{k-1}$. Hence assume in the rest that $t\le k$. 

If $t=1$, then $K(n,k)[T]$ is a disjoint union of $K_1$s and hence $|T|\le \binom{n-1}{k -1}$ by the Erd\H{o}s-Ko-Rado theorem. 

Suppose now $2\le t \le k$. Count the maximum number of $k$-subsets $A$, such that $A\cap B\ne \emptyset$ for any $B\in V(K)$. Since $A\cap B\ne \emptyset$ for any $B\in V(K)$, we can select one element from each $B\in V(K)$ and include it into $A$. Therefore we have  
$$\underbrace{\binom{k}{1}\binom{k}{1}\ldots \binom{k}{1}}_{t\mbox{-{\small times}}}\binom{n-t}{k-t}$$    
$k$-sets $A$, such that  $A\cap B\ne \emptyset$  for any $B\in V(K)$. Hence, using the inequality assume in the theorem, 
$$|T|\leq k^t\binom{n-t}{k-t}+t\leq \binom{n-1}{k-1}\,.$$ 
We conclude that $\gp(K(n,k))=\binom{n-1}{k-1}.$
\qed

Combining Theorem~\ref{thm:Kneser-special} with some further insights, the gp-number for case $k = 3$ can be determined.

\begin{theorem}  
\label{thm:kneser3}
If $n\ge 6$, then 
\[\gp(K(n,3)) = \left\{ 
{\begin{array}{{ll}}
20; & n = 6, \\
\binom{n-1}{2}; & n\geq 7\,.
\end{array}} \right.\]
\end{theorem}

\proof
Firstly, $K(6,3) = 10 K_2$, hence $\gp(K(6,3)) = 20$. In the rest we may thus assume that $n\ge 7$. Let $T$ be a general position set of $K(n,3)$. By Theorem~\ref{thm:gpsets}, every component of $K(n,3)[T]$ is a clique and let $Q$ be a largest such clique. 

If $n(Q)\ge 4$,  then $Q$ is the unique component of $T$ because no other vertex can have a non-empty intersection with all of the vertices of $Q$. 

Suppose $n(Q) = 3$ and note that $n\ge 9$ is necessary that this can happen. Then $K(n,3)[T]$ contains at most four components isomorphic to $K_3$. (This is indeed possible as demonstrated by blocks of a $2$-$(9,3,1)$ design.)
Suppose $K(n,3)[T]$ contains at least three components isomorphic to $K_3$. Then at most three vertices can have non-empty intersection with all the vertices in the $K_3$ components. Further, if there are exactly two components isomorphic to $K_3$, then each of the cliques allows $27$ further vertices to belong to $T$. The list of possible vertices intersects in only 6 vertices that can lie in $T$ besides the vertices from the two $K_3$.
Finally, suppose that there is only one component isomorphic to $K_3$, in 
which case again $27$ other vertices can belong to $T$. Every pair of disjoint vertices from this set of $27$ vertices excludes one vertex that has empty intersection with both sets. Therefore, at most $18$ vertices can lie in $T$ besides the vertices of the unique $K_3$, so at most $21$ in total. 

Assume next that $n(Q) = 2$. Then there are at most $10$ components isomorphic to $K_2$ in $T$, this is achieved by the $3$-subsets of $[6]$. If not all components are $K_2$, every other vertex in $T$ has to contain one element from each of the sets in $Q$. Then at most $3(n-2)$ other vertices are in $T$ (which we get by fixing one element from the first set, choosing one from the second set, and choosing the third one from remaining elements). Accordingly we have at most $\max \{20, 2+3(n-2)\}$ elements in the general position set of this shape. Finally, if all the components of $K(n,3)[T]$ are isomorphic to $K_1$, that is, if $T$ is an independent set, then $|T|\le \binom{n-1}{2}$ by the Erd\H{o}s-Ko-Rado theorem.

Consider now $K(7,3)$. Recall, that $\diam (K(7,3)) = 3$. A largest independent set in $K(7,3)$ has $15$ vertices, so $\gp(K(7,3)) \ge 15$. Suppose there is component isomorphic to $K_2$ in an arbitrary general position set $S$ with vertices $\{1,2,3\}$ and $\{4,5,6\}$. An addition vertex of $S$ must have non-empty intersection with $\{1,2,3\}$ and $\{4,5,6\}$. This vertex cannot have two elements in common with both vertices. If the intersections are not of the same cardinality, these three vertices will not be in general position. 
So an additional vertex must have one element in common with the vertices of $K_2$, 
from which we conclude that $\gp(K(7,3)) = 15$.

Next, let $n = 8$. Since $\diam(K(8,3)) = 2$, every independent set of $K(8,3)$ is a general position set, a largest such set has $\binom{7}{2} = 21$ vertices. Since $K(8,3)$ does not contain any subgraph isomorphic to $K_3$ and ${\rm max}(20,20) = 20$ (this $\max$ refers to the maximum two paragraphs above), we conclude that $\gp(K(8,3))=21$.

Let next $9\le n\le 19$. Then $K(n,3)$ contains an independent set of size $\binom{n-1}{2}$ and since $\diam(K(n,3)) = 2$, this is also general position set. Since for $n\ge 9$, $\binom{n-1}{2} \ge  2+3(n-2) > 20 > 18$, we conclude that $\gp(K(n,3)) = \binom{n-1}{2}$.

Finally, for every $n\ge 20$, the condition of Theorem~\ref{thm:kneser} is fulfilled, hence the assertion follows.\qed

\section{Cartesian products}
\label{sec:Cartesian}

In this section we prove a general lower bound on the gp-number of Cartesian product graphs. The bound is sharp as follows from the exact gp-number of the Cartesian product of two complete graphs. 

The {\em Cartesian product} $G\cp H$ of graphs $G$ and $H$ has the vertex set  $V(G)\times V(H)$ and the edge set $E(G\cp H)  = \{(g,h)(g',h'):\ gg'\in E(G)\mbox{ and } h=h', \mbox{ or, } g=g' \mbox{ and }  hh'\in E(H)\}$. If $(g,h)\in V(G\cp H)$, then the {\em $G$-layer}  $G^h$ through the vertex $(g,h)$ is the subgraph of $G\cp H$ induced by the vertices $\{(g',h):\ g'\in V(G)\}$. Similarly, the {\em $H$-layer} $^gH$ through $(g,h)$ is the subgraph of $G\cp H$ induced by the vertices $\{(g,h'):\ h'\in V(H)\}$. It is well-known that for given vertices $u = (g_1,h_1)$ and $v = (g_2,h_2)$ of $G\cp H$ we have ${\rm d}_{G\square H}(u,v) = {\rm d}_G (g_1,g_2) + {\rm d}_H (h_1,h_2)$. For more on the Cartesian product see the book~\cite{imrich-2008}.  

The announced lower bound reads as follows. 

\begin{theorem}
\label{thm:general-lower}
If $G$ and $H$ are connected graphs, then 
$$\gp(G \cp H) \geq \gp(G) + \gp(H) - 2\,.$$
\end{theorem}

\proof
Let $S_G\subseteq V(G)$ and $S_H\subseteq V(H)$ be gp-sets of $G$ and $H$,  respectively. Let $g\in S_G$ and $h \in S_H$. We claim that 
$$S = ((S_G \times \{h\}) \cup (\{g\} \times S_H)) \setminus \{(g,h)\}$$ 
is a general position set in $G\cp H$. 

Let $u,v\in S$. Suppose first that $u$ and $v$ lie in the layer $G^h$. 
Since layers in Cartesian products are convex, it follows that an arbitrary shortest $u,v$-path $P_{uv}$ lies completely in $G^h$. Since $G^h$ is isomorphic to $G$, it follows that $V(P_{uv}) \cap S = \{u,v\}$. Hence $(S_G \times \{h\}) \setminus \{(g,h)\}$ is a general position set in $G\cp H$. Analogously, $(\{g\} \times S_H) \setminus \{(g,h)\}$ is a general position set. 

Suppose now that $u = (g',h)\in G^h$, $v = (g,h')\in\, ^g\!H$, and let $P_{uv}$ be a shortest $u,v$-path in $G\cp H$. Suppose on the contrary that $P_{uv}$ contains some vertex $w$ of $S$ different from $u$ and $v$. We may without loss of generality assume that $w = (g'',h)$. Clearly, $g''\ne g'$. Furthermore, since $(g,h)\notin S$, we also have $g''\ne g$. Since the projection $P'$ of $P_{uv}$ on $G^h$ is a shortest path between $u = (g',h)$ and $(g,h)$ we infer that $P'$ passes through the vertex $(g'',h)$. This in turn implies that there exists a shortest $g',g$-path in $G$ that contains $g''$. This is a contradiction since $g$, $g'$, and $g''$ are pairwise different vertices. 

We have thus proved that $S$ is a general position set. Since $|S| = |S_G| + |S_H| - 2 = \gp(G) + \gp(H) - 2$ we are done. 
\qed

The bound of Theorem~\ref{thm:general-lower} is sharp as demonstrated by the equality case of the following result. 

\begin{theorem}
\label{thm:Hamming}
If $k\ge 2$ and $n_1, \ldots, n_k\ge 2$, then 
$$\gp(K_{n_1}\cp \cdots \cp K_{n_k}) \ge n_1 + \cdots + n_k - k\,.$$  
Moreover, $\gp(K_{n_1}\cp K_{n_2}) = n_1 + n_2 - 2$.  
\end{theorem}

\proof
To simplify the notation set $G = K_{n_1}\cp \cdots \cp K_{n_k}$. Let further $V(K_n) = [n]$, so that $V(G) = \{(j_1,\ldots,j_k):\ j_i\in [n_i], i\in [k]\}$.  

For $i\in [k]$ set $X_i = \{(1,\ldots, 1,j,1,\ldots, 1):\ j\in \{2,\ldots, n_i\}\}$, where $j$ is in the $i^{\rm th}$ coordinate. Clearly, $|X_i| = n_i - 1$. We claim that $X = \cup_{i\in [k]} X_i$ is a general position set of $G$.  

Let $u$, $v$, and $w$ be pairwise different vertices of $X$ and let $x\in X_p$, $v\in X_q$, and $w\in X_r$. If $p = q = r$, then $u$, $v$, and $w$ are in the same $K_{n_p}$-layer and thus induce a triangle. So they are in a general position. Suppose next that $p=q\ne r$. Then $d_G(u,v) = 1$, $d_G(u,w) = 2$, and $d_G(v,w) = 2$,  hence these three vertices are again in a general position in $G$. Finally, if $p|\ne q \ne r$, then $d_G(u,v) = d_G(u,w) = d_G(v,w) = 2$, and we have the same conclusion. This proves the claim. 

Since $X$ is a general position set and, clearly, $|X| = \sum_{i\in k}|X_i| = n_1 + \cdots + n_k - k$, the lower bound is proved. 

Let now $k=2$, so that $G = K_{n_1}\cp K_{n_2}$ and $V(G) = \{(i,j):\ i\in [n_1], j\in [n_2]\}$.  Since $\diam(G) = 2$, Theorem~\ref{thm:diameter2} applies. Clearly, $\omega(G) = \max \{n_1, n_2\}$. 

In the rest we are going to prove that $\eta(G) = n_1 + n_2 - 2$. We will prove this assertion by induction on $n_1+n_2$, the basic case $n_1 = n_2 =2$ being clear. Note also that if $n_2=2$ and $n_1\ge 3$, then the result also holds, that is, $\eta(G) = n_1$ in this case. 

Let $H$ be a complete multipartite subgraph of $\overline{G}$ and let $X_1, \ldots, X_k$ be the partite sets of $H$. We first claim that each $X_i$ is a subset of the vertex set of some layer. If $|X_1| = 1$ there is nothing to prove. Hence let $|X_1|\ge 2$ and suppose without loss of generality that $(1,1) \in X_1$. Since $X_1$ is an independent set, we have $(\{2,\ldots, n_1\} \times \{2,\ldots, n_2\}) \cap X_1 = \emptyset$. We may further suppose without loss of generality that $X_1$ contains another vertex from $K_{n_1}^{1}$, say $(i,1)$. Since $(i,1)$ is adjacent to all the vertices from $\{1\}\times \{2,\ldots, n_2\}$, we conclude that $X_1\subseteq V(K_{n_1}^{1})$. This proves the claim. 

By the above claim we may assume that $X_1 = \{(1,1), \ldots, (r,1)\}$, where $r\in [n_1]$. If $r = n_1$ then $k=1$ and hence $H$ has $n_1$ vertices. Since $n_2\ge 2$ we see that $n(H)\le n_1+n_2-2$. Suppose in the rest that $r<n_1$. Then none of the vertices from $(\{1,\ldots, r\} \times \{2,\ldots, n_2\})\cup (\{r+1,\ldots, n_1\} \times \{1\})$ lies in $H$. If follows that $X_2, \ldots, X_k$ lie in the subgraph induced by $\{r+1,\ldots, n_1\} \times \{2,\ldots, n_2\}$. The latter subgraph is isomorphic to $K_{n_1-r}\cp K_{n_2-1}$. If $n_1-r\ge 2$ and $n_2-1\ge 2$, then by the induction (on $n_1+n_2$) we have  $\eta(K_{n_1-r}\cp K_{n_2-1}) = (n_1 - r) + (n_2 -1) - 2 = n_1 + n_2 - r - 3$. It follows that $n(H) = n_1 + n_2 - 3$. Let now $n_1-r\le 1$, which in turn implies by the above that $n_1-r=1$, that is $r=n_1-1$. But then $X_1 = \{(1,1), \ldots, (n_1-1,i)\} \cup \{(n_1,2), \ldots, (n_1,n_2)\}$ induce a complete bipartite graph in $\overline{G}$, hence $\eta(G) \ge (n_1-1) + (n_2 - 1) = n_1 + n_2 - 2$. Suppose finally that $n_2-1\le 1$, that is, $n_2\le 2$, and so $n_2=2$, the case that was already considered. 
\qed

Note that the lower bound of Theorem~\ref{thm:Hamming} for at least three factors is stronger than the bound one can deduce by induction from Theorem~\ref{thm:general-lower}.

\section{The gp-number of some graph operations}
\label{sec:more-opearations}

In this section we consider the gp-number of joins of graphs, of coronas over graphs, and of line graphs. For this sake the following concept will be useful. Complete subgraphs $Q$ and $Q'$ in a graph $G$ are {\em independent} if $d_G(u,u')\ge 2$ for every $u\in V(Q)$ and every $u'\in V(Q')$. (This concept has been very recently introduced and applied in~\cite{chadha-2018+}.) Note that the complete subraphs from Theorem~\ref{thm:gpsets} are independent by definition. Setting $\rho(G)$ to denote the maximum number of vertices in a union of pairwise independent complete subgraphs of $G$, we have: 

\begin{theorem}
\label{thm:diameter-2-second}
If $\diam(G) = 2$, then $\gp(G) = \rho(G)$. 
\end{theorem}

\proof
Let $G$ be a graph of diameter $2$. Clearly, $\rho(G) \ge \omega(G)$ and $\rho(G)\ge \eta(G)$. Theorem~\ref{thm:diameter2} thus implies that $\rho(G) \ge \gp(G)$. Conversely, the $\rho(G)$ vertices from a largest union of pairwise independent cliques form a general position set by Theorem~\ref{thm:gpsets}. Therefore, $\gp(G)\ge \rho(G)$.   
\qed

The reason that in Theorem~\ref{thm:diameter-2-second} $\gp(G)$ is expressed only with $\rho(G)$, while in Theorem~\ref{thm:diameter2} two invariants are used, is that $\rho(G)$ encapsulates $\omega(G)$ while $\eta(G)$ does not. 

\subsection{Joins and coronas}

If $G$ and $H$ are disjoint graphs, then the {\em join} $G+H$ of $G$ and $H$ is the  graph with the vertex set $V(G+H) = V(G)\cup V(H)$ and the edge set $E(G + H) = E(G)\cup E(H)\cup \{xy:\ x\in V(G), y\in V(H)\}$. If both $G$ and $H$ are complete, so it is $G+H$ and hence $\gp(G + H) = \gp(K_{n(G)} + K_{n(H)}) = \gp(K_{n(G) + n(H)}) = n(G + H)$. Otherwise, that is, if at least one of $G$ and $H$ is not complete, then $\diam(G+H) = 2$. In this case we have: 

\begin{proposition}
\label{prop:join}
If $G$ and $H$ are graphs, then 
\begin{eqnarray*}
\gp(G + H ) & = & \max\{\omega(G) + \omega(H), \eta(G), \eta(H)\}\\
& = & \max\{\omega(G) + \omega(H), \rho(G), \rho(H)\}\,.
\end{eqnarray*}
\end{proposition}

\proof
Since  $\diam(G+H) = 2$, Theorem~\ref{thm:diameter2} applies. It is straightforward that $\omega(G + H) = \omega(G) + \omega(H)$ and that $\eta(G + H) = \max \{\eta(G), \eta(H) \}$. Hence the first equality. 

A complete subgraph $Q$ of $G + H$ lies completely in $G$, or completely in $H$ or is a join of a complete subgraph of $G$ and a complete subgraph of $H$. If $Q$ is of the latter form, then it is at distance $1$ to every other complete subgraph of $G+H$. If follows that $\rho(G+H) = \max\{\omega(G) + \omega(H), \rho(G), \rho(H)\}$. The second equality then follows by Theorem~\ref{thm:diameter-2-second}. 
\qed

Let $G$ and $H$ be graphs where $V(G) = \{v_1, \ldots ,v_{n(G)}\}$. The {\em corona} $G\circ H$ of graphs $G$ and $H$ is obtained from the disjoint union of $G$ and $n(G)$ disjoint copies of $H$, say $H_1,\ldots, H_{n(G)}$, where for all $i\in [n(G)]$, the vertex $v_i\in V(G)$ is adjacent to each vertex of $H_i$. 

\begin{theorem}
\label{thm:corona}
If $G$ and $H$ are graphs where $n(G)\ge 2$, then 
$$\gp(G\circ H) = n(G)\rho(H)\,.$$ 
\end{theorem}

\proof
Let $V(G) = \{v_1, \ldots ,v_{n(G)}\}$ and let  $H_1,\ldots, H_{n(G)}$ be the corresponding copies of $H$ in $G\circ H$. Note first that the statement is clear for the corona $K_2\circ K_1 = P_4$. So we may assume in the rest that if $n(G) = 2$ then $n(H)\ge 2$. 

Let $S$ be a gp-set of $G\circ H$. Suppose first that $S\cap V(G)\ne \emptyset$. We may assume without loss of generality that $v_1\in S$. If there exists a vertex $w\in S\cap V(H_1)$, $w\ne v_1$, then for any vertex $x\in V(G\circ H) \setminus (V(H_1)\cup\{v_1\})$, the vertex $v_1$ lies on a shortest $w,x$-path. Consequently, $S\subseteq V(H_1)\cup\{v_1\}$. If $n(G) = 2$, then since we have assumed $n(H)\ge 2$, the union of a gp-set of $H_1$ and a gp-set of $H_2$ has cardinality bigger that $S$ because $\gp(H)\ge 2$. And if $n(G)\ge 3$, we get a similar contradiction. It follows that if $v_1\in S$, then $S\cap V(H_1) = \emptyset$. But then $S' = S \cup \{w\}\setminus \{x_1\}$, where $w$, is an arbitrary vertex of $H_1$ is also a gp-set. In summary, we have proved that we may without loss of generality assume that $S\cap V(G) = \emptyset$.   

So let now $S$ be a gp-set of $G\circ H$ with $S\cap V(G) = \emptyset$. By Theorem~\ref{thm:gpsets}, the components of $(G\circ H)[S]$ are independent complete graphs. Hence, $S$ restricted to $H_i$ has at most $\rho(H)$ vertices. On the other hand, since independent complete subgraphs of $H_i$ are pairwise at distance $2$, they form (in view of Theorem~\ref{thm:gpsets}) a general position set. But then taking such complete subgraphs in every $H_i$ yields a general position set of order $n(G)\rho(H)$. 
\qed

\subsection{Line graphs of complete graphs}

If $G$ is a graph, then the {\em line graph} $L(G)$ of $G$ is the graph with $V(L(G)) = E(G)$, two different vertices of $L(G)$ being adjacent if the corresponding edges share a vertex in $G$. 

\begin{theorem} 
If $n\ge 3$, then 
\[\gp(L(K_n)) = \left\{ {\begin{array}{ll}
n; & 3\mid n\,, \\
n-1; & 3\nmid n\,.
\end{array}} \right.\]
\end{theorem}

\proof
Let $n\ge 3$ and $V(K_n) = [n]$. To simplify the notation set $G_n = L(K_n)$. Since $\omega(G_n)=n-1$, we have $\gp(G_n) \geq n-1$. 

We next claim that $\gp(T(n))\leq n$. Let $S$ be a $\gp$-set of $G_n$  and let $K_{n_1}, \ldots, K_{n_k}$ be the connected components of $G_n[S]$, so that $\gp(G_n) = |S| = n_1 + \cdots + n_k$. A vertex $u$ of $G_n$ corresponds to an edge of $K_n$, that is, to a pair of vertices $\{j,j'\}$ and we may identify $u$ with $\{j,j'\}$. Using this convention, for $i\in [k]$ set 
$$X_i = \bigcup_{\{j,j'\}\in V(K_{n_i})} \{j, j'\}\,.$$ 
Since the complete subgraphs $K_{n_i}$ are pairwise independent, it follows that if $i\ne i'$, then $X_i\cap X_{i'}=\emptyset$.  Setting $x_i=|X_i|$ we infer that $x_i\geq n_i$ and hence 
\begin{equation}  
\label{eq:n_i-x_i}
\gp(G_n) = |S| = n_1 + \cdots + n_k \leq x_1 + \cdots + x_k \leq n\,,
\end{equation}
and the claim is proved.

If $3\mid n$, then 
$$S=\{\{{3i+1}, {3i+2}\}, \{{3i+1}, {3i+3}\}, \{{3i+2}, {3i+3}\}:\ 0\leq i \leq \frac{n}{3}-1\}$$
is a gp-set of $G_n$, and hence  $\gp(G_n)= n$. 

Suppose now that $3\nmid n$. Then at least one $n_i \ne 3$ and for it we have $n_i < x_i$. In view of~\eqref{eq:n_i-x_i} this means that $\gp(G_n) < n$. As we have already observed that  $\gp(G_n) \geq n-1$, the argument is complete. 
\qed

\section*{Acknowledgements}

We acknowledge the financial support from the Slovenian Research Agency 
(research core funding No.\ P1-0297 and projects J1-9109, N1-0095).


\end{document}